 \documentclass{ajour}
\input xy
\xyoption{all}

\usepackage{amsmath}
\usepackage{latexsym}
\usepackage{amssymb} 
\usepackage{amscd} 
\usepackage{psfig} 
\usepackage{multicol}

\def\P{\Bbb{P}} 
\def\R{\Bbb{R}} 
\def\C{\Bbb{C}} 

\def\iso{\cong}

\begin{document}

%


\authorrunninghead{R. F. Goldin}
 \titlerunninghead{Weight Varieties}





\title{The Cohomology Ring of Weight Varieties and Polygon
Spaces}


\author{R. F. Goldin\thanks{Supported by a National Science Foundation
Postdoctoral Fellowship}\thanks{The author would like to
thank A. Knutson and S. Billey for many useful conversations about
this work, and the referee for such detailed and helpful comments.}}
\affil{University of Maryland}

\email{goldin@math.umd.edu}

%







\abstract{We use a theorem of Tolman and Weitsman
\cite{TW:abelianquotient} to find explicit formul\ae\ for the rational
cohomology 
rings of the symplectic reduction of flag varieties in $\C^n$, or
generic coadjoint orbits of 
$SU(n)$,  by (maximal) torus actions. We also calculate the cohomology
ring of the moduli space of $n$ points in $\C P^k$, which is
isomorphic to the Grassmannian of $k$ planes in $\C^n$, by realizing
it as a degenerate coadjoint orbit.}


\begin{article}


\section{Introduction}\label{se:intro}
 For $M$ a manifold with a Hamiltonian $T$ action and moment map
$\phi:M\rightarrow\mathfrak{t}^*$, the symplectic reduction is defined
as $$M/\!/T(\mu):=\phi^{-1}(\mu)/T$$ for any regular value $\mu$ of
$\phi$. {\em Weight varieties} are a special case of symplectic
reduction. Let $M=\mathcal{O}_\lambda$, a coadjoint orbit of a
compact, semi-simple Lie group $G$ through the point
$\lambda\in\mathfrak{t}^*$, and consider the action of the maximal
torus $T\subset G$ on $\mathcal{O}_\lambda$. If $G=SU(n)$, we identify the set of Hermitian matrices $\mathcal{H}$ with 
$\mathfrak{g}^*$ by $tr:A \rightarrow
i\cdot Trace(A\cdot)$ for all $A\in \mathcal{H}$. Under this identification, we can think of
$\lambda$ as a matrix with real diagonal entries
$(\lambda_1,\dots,\lambda_n)$, and 
 $\mathcal{O}_\lambda$ as an adjoint orbit of $G$ through $\lambda$.
The moment map for the $T$ action on $\mathcal{O}_\lambda$
takes a matrix to its diagonal entries. Thus
$\mathcal{O}_\lambda/\!/T(\mu)$ consists of Hermitian matrices
with spectrum $\lambda$ and
diagonal entries $\mu$, quotiented out by the action of
diagonal matrices. The symplectic reduction
$\mathcal{O}_\lambda/\!/T(\mu)$ is a weight variety.

The generic coadjoint orbit of $SU(n)$ is symplectomorphic to the
complete flag variety in $\C^n$ with a symplectic structure given by
the spectrum of the orbit. Degenerate coadjoint orbits are homeomorphic
to flag varieties with various dimensions missing: for example if the
spectrum consists of only two values, then the coadjoint orbit is
homeomorphic to $Gr(k,n)$, the Grassmannian of $k$ planes in $\C^n$.
Then $\mathcal{O}_\lambda/\!/T(\mu)$ is the symplectic reduction of a
flag variety by a torus; however, the reduction depends on the
symplectic structure on the flag variety.

Weight varieties have appeared in several different contexts. They
were first termed as such by Knutson in \cite{Knutson:thesis} because
of their relationship to the weight spaces of
representations. Irreducible representations $V$ of complex $G$ are
realized as the holomorphic sections of line bundles over the flag
varieties. The dimension of the {\em weight spaces}, or the
irredicible representations of $T$ in $V$ (specified by the weight of
the $T$ action on each isotypic component of $V$) is the {\em quantization} of
the symplectic reduction of flag varieties \cite{GS:quantization};
hence these reductions were named weight varieties.  Techniques to
compute their Betti numbers have been developed by Kirwan
\cite{KiBook:quotients} and Klyachko \cite{Klyachko:Bettinumbers}. The
advantage of the results in this article is that one obtains the
cohomology ring structure. For Betti numbers, one may find that the
calculations are straightforward (and hence programmable) using
Theorems \ref{th:generic} and \ref{th:grassmannian}, but other methods
are likely more efficient.  In a small number of cases, the spaces
have been found explicitly \cite{GLS:fibrations}.

The symplectic reduction of the Grassmannian of $k$-planes in $\C^n$
also arises as the moduli space of $n$ points
on $\C P^k$, as we see below; more details for the case when $k=1$ can
be found in the work 
of Klyachko \cite{Klyachko:Bettinumbers}, Hausmann and Knutson \cite{Knutson:polygonspaces} and Kapovich and Millson
\cite{Millson:euclideanpolygons}. The integer cohomology ring of the $k=1$
polygon spaces was computed in 
\cite{Knutson:cohomologypolygonspaces}.

The generic case we consider is the coadjoint orbit of $SU(n)$ of
matrices with a specified set of distinct eigenvalues. Equivalently,
it is the orbit through
$\lambda\in\mathfrak{t}^* := Lie(T)^*$,
where $\lambda$ consists of $n$ distinct eigenvalues
$\lambda=(\lambda_1,\dots,\lambda_n)$  with $\sum_1^n\lambda_i=0$, and
$T$ is the $(n-1)$-dimensional maximal torus of $SU(n)$. 
We order them such that
$\lambda_1>\dots >\lambda_n$. These coadjoint orbits are diffeomorphic
to the complete flag manifold $Fl(n)$ in $\C^n$.
The Grassmannian of $k$-planes in $\C^n$ is diffeomorphic to the $SU(n)$
coadjoint orbit through $\nu\in\mathfrak{t}^*$ where
$\nu=(\nu_1,\dots,\nu_1,\nu_2,\dots,\nu_2)$ consists of two distinct
eigenvalues $\nu_1$ and $\nu_2$. The dimension of the $\nu_1$
eigenspace is $k$ and that of the $\nu_2$ eigenspace is $n-k$, so we
have $k\nu_1+(n-k)\nu_2 =0$. We will write $Gr(k,n)_\nu$ to indicate
these degenerate coadjoint orbits.

We use the notation $\lambda_w$ to indicate the point $w\lambda
w^{-1}\in \mathfrak{t}^*$
for any permutation
$w$ on $n$ letters. This has the unfortunate consequence that
$\lambda_w = (\lambda_{w^{-1}(1)}, \dots, \lambda_{w^{-1}(n)})$ but allows
us to use left actions consistently throughout the paper. Furthermore,
let $\Delta(x,u)\in\C[x_1,\dots,x_n,u_1,\dots,u_n]$ be the polynomial
$$
\Delta (x,u) = \prod_{i<j}(x_i-u_j).
$$ $\Delta(x,u)$ is sometimes called the {\em
determinant polynomial}. 

Before we state the main theorems, we briefly introduce {\em divided
difference operators}, which 
are explained in detail in Section \ref{se:schubertpolys}. 
\begin{definition}
Let $f(x_i,x_{i+1})$ be a polynomial of variables $x_i$ and $x_{i+1}$
and possibly other variables. For each $1\leq i\leq n$ the
{\em divided difference operator $\partial_i$ associated to $i$} acts on $f$ as follows:
$$
\partial_i f (x_i,x_{i+1})  = \frac{f(x_i,x_{i+1})-f(x_{i+1},x_i)}{x_i-x_{i+1}}.
$$
\end{definition}
These simple divided difference operators take polynomials of degree $k$ to
polynomials of degree $k-1$. For any element $w\in S_n$, the
permutation group on $n$ letters, one can associate a divided
difference operator $\partial_w$ as follows. Write $w$ as a product of
simple transpositions $s_{i_1}\dots s_{i_l}$ where $s_{i_j}$ is the
simple trasposition that switches $i_j$ and $i_j+1$. For any such
product where $l$ is minimal, the composition
$$\partial_w := \partial_{i_1}\cdots\partial_{i_l}$$
is well-defined (see \cite{MacDonald}).
We are now ready to state the main theorems of this article.

\begin{theorem}\label{th:generic}
Let $\mathcal{O}_\lambda$ be a generic coadjoint orbit of $SU(n)$.
The rational cohomology of $\mathcal{O}_\lambda/\!/T(\mu)$ is
isomorphic to the ring
$$
\frac{\C [x_1,\dots,x_n,u_1,\dots, u_n]}{\big(\prod_{i=1}^n(1+u_i)-\prod_{i=1}^n(1+x_i), \sum_{i=1}^nu_i, \partial_v\Delta(x,u_\tau)\big)}  
$$ 
for all $v,\tau \in S_n$ such that $\sum_{i=k+1}^{n} \lambda_{v(i)} <
 \sum_{i=k+1}^n\mu_{\tau(i)}$ for some $k=1,\dots,n-1$. Here
 $\deg x_i=\deg u_i =2$, and $\prod(1+u_i)-\prod(1+x_i)$ is the graded
 difference of symmetric functions in the $x_i$s and $u_i$s.
\end{theorem}

There is a similar statement for the degenerate case of the
Grassmannian, in which there are only two distinct eigenvalues.

\begin{theorem}\label{th:grassmannian}
The rational cohomology of $Gr(k,n)_\nu/\!/T(\mu)$ is isomorphic to the ring
$$
\frac{\C[\sigma_i(x_1,\dots, x_k), \sigma_i(x_{k+1},\dots, x_n), u_1,\dots,
u_n]}{\big(\sigma_i(x_1,\dots, x_n)    
-\sigma_i(u_1,\dots, u_n), \sum_{i=1}^nu_i, \partial_v\Delta(x,u_\tau)\big)} 
$$
for all $v, \tau$ such that $\sum_{i=k+1}^n \lambda_{v(i)}<\sum_{i=k+1}^n 
\mu_{\tau(i)}$ for some $k$, and $\partial_v\Delta(x,u_\tau)$ is symmetric in
$x_1,\dots, x_k$ and in $x_{k+1},\dots, x_n$.
Here $\sigma_i$ are
the symmetric polynomials in the indicated variables and $\deg
x_i=\deg u_i =2$. Note that $\sigma_i(x_1,\dots, x_n)    
-\sigma_i(u_1,\dots, u_n)$ for all $i$ is equivalent to
$\prod_{i=1}^n(1+u_i)-\prod_{i=1}^n(1+x_i)$.
\end{theorem}

The forgetful map $Fl(n)\rightarrow Gr(k,n)$ which ``remembers''
only the $k$-planes of each flag is a $T$-equivariant map which
induces an injection on equivariant cohomology:
\begin{equation}\label{eq:grasstoflag}
H_T^*(Gr(k,n))\hookrightarrow H_T^*(Fl(n)).
\end{equation}
The classes in $H_T^*(Fl(n))$ which are symmetric in $x_1,\dots, x_k$
and in $x_{k+1},\dots, x_n$ are precisely those 
in the image of the map (\ref{eq:grasstoflag}).

Theorem \ref{th:grassmannian} allows one to compute the cohomology of
the moduli 
space of $n$ points in $\C P^{k-1}$, as the reduction of the Grassmannian
is isomorphic to this moduli space \cite{GM:geometry},
\cite{GGMS:combinatorial}. In \cite{Knutson:polygonspaces}, Knutson
and Hausmann observe that this Gelfand-MacPherson correspondence is
just a dual pair symplectic reduction. 
The Grassmannian $Gr(k,n)$ is
realized as a symplectic reduction of $\C^{nk}$ by a $U(k)$ action on
the right. Then by reducing in stages we have
\begin{align*}
Gr(k,n)/\!/T^n &=(\C^{nk}/\!/U(k))/\!/T^n= \C^{nk}/\!/(T^n\times U(k))\\
&=(\C^{nk}/\!/T^n)/\!/U(k)= \prod_{i=1}^{n}\C P^{k-1}/\!/U(k).
\end{align*}
This last space is exactly the moduli space of $n$ points in $\C P^k$.
The group $U(k)$ does not act effectively on $\prod_{i=1}^{n}\C
P^{k-1}$. The center, consisting of scalar matrices, acts
trivially. We quotient this $S^1$ out and find
$$
Gr(k,n)/\!/T = \prod_{i=1}^{n}\C P^{k-1}/\!/PU(k),
$$
where $T$ is the $(n-1)$-dimensional torus used in Theorem
\ref{th:grassmannian}. 
Putting in the symplectic structure, we make the
following statement.
\begin{corollary} Let $\mathcal{M}$ be the moduli space of $n$ points
in $\C P^{k-1}$, realized as above by a symplectic reduction as follows:
\begin{align*}
\mathcal{M} &=  \big(\C
^{nk}/\!/U(k)\big)\big(aI\big)/\!/T^n(\mu_1+a,
\mu_2+a, \dots, \mu_n+a)\\
&= Gr(k,n)_{(\nu_1,\dots,\nu_1,\nu_2,\dots, \nu_2)}/\!/T(\mu_1,\dots, \mu_n).
\end{align*}
where $a=(\nu_1-\nu_2)$ and $I$ is the identity matrix in $\mathfrak{u}^*(n).$
Then $H^*(\mathcal{M})$ is isomorphic to the ring
$$
\frac{\C[\sigma_i(x_1,\dots, x_k), \sigma_i(x_{k+1},\dots, x_n), u_1,\dots,
u_n]}{\big(\sigma_i(x_1,\dots, x_n)    
-\sigma_i(u_1,\dots, u_n), \sum_{i=1}^nu_i, \partial_v\Delta(x,u_\tau)\big)} 
$$
for all $v, \tau$ such that $\sum_{i=k+1}^n \lambda_{v(i)}<\sum_{i=k+1}^n 
\mu_{\tau(i)}$ for some $k$, and $\partial_v\Delta(x,u_\tau)$ is symmetric in
$x_1,\dots, x_k$ and in $x_{k+1},\dots, x_n$.
\end{corollary}

\begin{remark}
The equivariant Chern class of the tangent bundle $TM\rightarrow M$
descends to the ordinary (total) Chern class of $M/\!/T$ under symplectic
reduction by a torus. We find the Chern class of the tangent bundle to 
$Fl(n)\iso
\mathcal{O}_\lambda$. We note that $c(\C P^n) = (1-x)^{n+1}$, where
$x$ is the first Chern class of the tautoligical line bundle $S$ over
$CP^n$ (see \cite{BottTu}).
Using an inductive argument on the fibration
$Fl(n-1)\hookrightarrow Fl(n)\rightarrow \C P^n$, one can show that in
the basis used above, 
the Chern class of $\mathcal{O}_\lambda$ is  $(1-x_1)^n(1-x_2)^{n-1}\cdots
(1-x_{n-1})^2$. Then this is also the total Chern class of the weight
varieties $\mathcal{O}_\lambda/\!/T(\mu)$.
\end{remark}

There are two essential facts that come into play in the results presented
here. For
$M$ a symplectic manifold with a Hamiltonian torus action, there is a
restriction map in equivariant cohomology from $M$
to the $\mu$-level set $\phi^{-1}(\mu)$ of the moment map $\phi$. The
(rational) equivariant cohomology of the level set is equal to the
regular cohomology of the reduced space $M/\!/T(\mu):=
\phi^{-1}(\mu)/T$. The first theorem is that the resulting map is a
surjection.

\begin{theorem}[Kirwan]\label{th:surjection}
Let $M$ be a Hamiltonian $T$ space with moment map $\phi$ and
$\mu\in\mathfrak{t}^*$ a regular value of $\phi$. Then the map induced by
restriction to the level set $\phi^{-1}(\mu)$
$$
\kappa_\mu: H_T^*(M)\longrightarrow H^*(M/\!/T(\mu))
$$ is a surjection.
\end{theorem}

Secondly, the restriction map in equivariant cohomology induced by
the inclusion of the fixed point set $M^T$ into $M$
is  an injection.

\begin{theorem}[Chang-Skjelbred, Kirwan]\label{th:injection}
Let $M$ be a Hamiltonian $T$ space with fixed point set $M^T$. The
natural map
$$
r^*:H_T^*(M)\hookrightarrow H_T^*(M^T)
$$
is an inclusion. 
\end{theorem}
Note how different this is from ordinary cohomology. If $M$ has
isolated fixed points, for example, $H^*(M^T)$ is zero except in
degree 0, yet $H^*(M)$ may have cohomology in higher degree.

Theorem \ref{th:injection} suggests the following definition.
\begin{definition}
Let $\alpha\in H_T^*(M)$ be an equivariant cohomology class on a
compact Hamiltonian $T$ space with fixed point set $M^T$. Define the {\em support of $\alpha$} to
be $$ supp\ \alpha = \{C \mbox{\em \ connected\ component\ of\ } M^T| r_C^*(\alpha)\neq 0\}$$
where $r_C^*:H_T^*(M)\rightarrow H_T^*(C)$ is the restriction to the
equivariant cohomology of the fixed component $C$.
\end{definition}

By Theorem \ref{th:surjection}, the cohomology of the symplectic
reduction can be computed as the quotient of the equivariant
cohomology $H_T^*(M)$ by the kernel of the Kirwan map
$\kappa$. Theorem \ref{th:injection} indicates that the kernel may be
generated by cohomology classes which have certain properties
restricted to the fixed point set. This line of reasoning was
exploited by Tolman and Weitsman who described the kernel $\kappa_\mu$
\cite{TW:abelianquotient}. 

\begin{theorem}[Tolman-Weitsman]\label{th:tolweis}
Let $M$ be a compact symplectic manifold with a Hamiltonian $T$
action. Let $\phi:M\longrightarrow
\mathfrak{t}^*$ be a moment map such  
that $\mu$ is a regular value. Define
$$
M^\mu_\xi:=\{m\in M|\langle \phi(m),\xi\rangle \leq \langle\mu,\xi\rangle)\}
$$
and
$$
K_\xi:=\{\alpha\in H^*(M)| supp\ \alpha \subset M^\mu_\xi\}.
$$
Then the kernel of the natural map $\kappa_\mu:H_T^*(M)\longrightarrow
H^*(M/\!/T(\mu))$ is the ideal $\langle K\rangle$ generated by 
$$
K:=\bigcup_\xi K_\xi.
$$
\end{theorem}
The classes in the kernel are generated by classes
which have non-zero restriction to fixed points which (under the
moment map) lie entirely to one side of a hyperplane $\xi^\perp_\mu$
through $\mu$ in $\mathfrak{t}^*$.  Because $M$ is compact, 
only a finite number of hyperplanes will be necessary.

The contribution of this article is the application of the work of
Tolman and Weitsman to the case where $M$ is a coadjoint orbit of $SU(n)$. Here there is an explicit description of a
generating set of classes for the ideal $\langle K\rangle$, which
allows one to actually compute the cohomology ring of weight
varieties. These classes are represented by double Schubert
polynomials,
permuted by the Weyl group and  
satisfying certain properties. Double Schubert polynomials were first introduced by Lascoux and Sch\"utzenberger
\cite{LS:Schubertpolynomials}, \cite{MacDonald}.  As a corollary, for
$SU(n)$ coadjoint orbits, the only vectors $\xi\in\mathfrak{t}$ needed for
Theorem \ref{th:tolweis} are fundamental weights and their
permutations (Corollary \ref{co:sufficient}). Equivalently, the kernel of the Kirwan map $\kappa_\mu$
is generated by classes which restrict to zero on one side of
hyperplanes $\xi^\perp := ann(\xi)$, translated to contain $\mu$,
parallel to codimension-one walls of the moment polytope.

\section{The $T$-Equivariant Cohomology of Flag Manifolds and
Grassmannians}\label{se:flag} 

There are two descriptions of the equivariant cohomology of the flag
manifold that we will use here. First is a presentation of the ring as
a quotient of a polynomial ring by relations.  This is due in the
nonequivariant setting to Borel \cite{Borel:equivariant} and can
be found in \cite{BottTu}. The equivariant version can be obtained
from the standard computation by using the Borel construction: By definition,
$H_T^*(Fl(\C^n)):= H^*(Fl(\C^n)\times_TET)$ where $ET$ is a universal
$T$ bundle, and note that $Fl(\C^n)\times_TET$ is a flag bundle over
the classifying space $BT:=ET/T$. Using standard methods, one computes
the cohomology of this bundle. Details can be found in
\cite{Goldin:thesis}, \cite{Brion:equivariantcohomologynotes}.

The second description is due to Arabia \cite{Arabia:equivcohom}, and
says that there is a basis for $H_T^*(G/T)$ as a module over
$H_T^* := H_T^*(pt)$ which has certain properties that we expand upon
below. The connection between these two pictures in the case that $G=SU(n)$ is provided by double
Schubert polynomials. We show by Theorem \ref{th:corresp} and Theorem
\ref{th:integral} that double Schubert polynomials, and their
corresponding classes under a quotient map, have the properties
specified by Arabia. Hence the linear basis can be expressed in the
presentation due to Bernstein, Gelfand, and Gelfand.

\subsection{Presentations of $H_T^*(Fl(\C^n))$ and
$H_T^*(Gr(k,n))$}\label{se:flagpresentation} 
The $T^n$ action on $Fl(\C^n)$ which is induced by the action on
$\C^n$ of weight 1 on each of $n$ copies of $\C$ is not effective: a
diagonal $S^1\subset T^n$ fixes $Fl(\C^n)$. We quotient by this circle
and let $T$ be the $(n-1)$-torus which acts effectively. 
\begin{theorem} \label{th:cohFlag}
Let $T$ be the maximal $(n-1)$-dimensional torus in $SU(n)$.
The $T$-equivariant cohomology of the generic $SU(n)$-coadjoint orbit is 
$$
H_T^*(Fl(\C^n))=\frac{\C[x_1,\dots ,x_n,u_1,\dots,u_n]}
{\big(\prod_{i=1}^n(1+u_i)-\prod_{i=1}^n(1+x_i),\
\sum_{i=1}^nu_i\big)}, \ \ \ i=1,\dots,n
$$
where  $\deg u_i=\deg x_i=2$, and the $u_i$ are the image of the
classes from the module structure $H_T^*\rightarrow H_T^*(Fl(\C^n)).$
\end{theorem}
It is straightfoward to see that $H_T^*(Fl(\C^n)) = H^*(Fl(E))$, the
flags of the bundle $E:= \C^n\times_TET\rightarrow BT$ for a certain
action of $T$ on $\C^n$. $Fl(E)$ can be canonically realized as a tower
of projective bundles over $BT$
$$
\xymatrix{
Fl(E)=\P (Q_{n-1}) \ar[r]^{\pi_{n-1}} & \dots \ar[r]^{\pi_3} & \P
(Q_1) \ar[r] ^{\pi_2} 
& \P(E) \ar[r]^{\pi_1} & BT
}
$$
where $Q_1:= \pi_1^*E/S_1$ for $S_1$ the tautological line bundle over
$\P(E)$. Thus $\pi_1^*E = Q_1\oplus S_1$ and the bundle
$Q_1\rightarrow \P (E)$ can to be projectivized to form $\pi_2:\P
(Q_1) \rightarrow \P (E)$. We repeat this process to obtain $Q_k:= \pi_k^*\dots\pi_1^*E/S_1\oplus \dots\oplus S_k$
where $S_k$ is the tautological line bundle over $\P(Q_{k-1})$. 
The pullback of $E$ to $Fl(E)$ canonically splits into a sum of line
bundles
\begin{equation*}
\begin{CD}
\pi^*E = S_1\oplus \dots\oplus S_n @>>> E\\
@VVV @VVV\\
Fl(E)@>\pi>> BT
\end{CD}
\end{equation*}
By definition, $x_i= c_1(S_i)$.
It follows that
$$ H_T^*(Fl(E)) = H^*(BT)[x_1,\dots, x_n]/\prod_{i=1}^n(1+x_i) = c(E).
$$
Using the splitting principle one obtains that the total Chern class
of $E$ is $\prod_{i=1}^n (1+u_i)$, where $u_i$ is the first Chern
class of the tautological line bundle over the $i$th copy of $S^1$ in
$T$ under a choice of decomposition $T=S^1\times\cdots \times
S^1$. Then restricting to the $n-1$ dimensional torus that acts
effectively on $Fl(E)$, we get the desired result.

\begin{theorem}\label{th:cohGrass}
Let $T$ be the maximal $(n-1)$-dimensional torus in $SU(n)$. The
$T$-equivariant cohomology of the Grassmannian of complex $k$-planes
in $\C^n$ is a subring of
$H_T^*(Fl(\C^n))$ and can be written
$$
H^*_T(Gr(k,n)) = \frac{\C[\sigma_i(x_1,\dots ,x_k), \sigma_i(x_{k+1},\dots ,x_n),u_1,\dots,u_n]}
{\big(\sigma_i(x_1,\dots, x_n)-\sigma_i(u_1,\dots, u_n),\ \sum_{i=1}^nu_i\big)},
$$
where $\sigma_i$ is the $i^{th}$ symmetric
function in the indicated variables, $i=1,\dots, n$. The Chern classes
$x_i$ and $u_i$ are the same as those in Theorem $\ref{th:cohFlag}$.
\end{theorem}

We now explore the relationship between this quotient description of
the equivariant cohomology, and the description as a subring of the
equivariant cohomology of the fixed point set.

\subsection{The Restriction of Equivariant Cohomology to Fixed Points}\label{se:restriction}
By Theorem \ref{th:injection}, the $T$-equivariant cohomology of
$Fl(\C^n)$ can be embedded in a direct sum of polynomial rings. We
calculate this restriction explicitly using geometric means. A purely
algebraic proof can be found in
\cite{Brion:equivariantcohomologynotes}. In
\cite{Billey:kostant} Billey finds a simple formula for the restriction of
Kostant polynomials (generalized double Schubert polynomials) to the fixed point set. 

The fixed point set $Fl(\C^n)^T$ is indexed by the Weyl group
$W=N(T)/T$. The $T$ action on $V=\C^n$ splits $V$ into a 
sum of 1-dimensional vector spaces, or lines which we order and call $l_1,\dots,l_n$. The 
fixed points of $T$ on $Fl(V)$ are the flags which can be written
$\langle l_{i_1}\rangle\subset\langle 
l_{i_1},l_{i_2}\rangle\subset\dots\subset\langle 
l_{i_1},\dots,l_{i_n}\rangle=V$. We call
$\langle l_{1}\rangle\subset\langle 
l_{1},l_{2}\rangle\subset\dots\subset\langle 
l_{1},\dots,l_{n}\rangle=V$ the base flag and label the
fixed points by the corresponding permutation in $W$:
$$p_w:=\langle 
l_{w(1)}\rangle\subset\langle l_{w(1)}, 
l_{w(2)}\rangle\subset \cdots\subset\langle l_{w(1)}, \dots, 
l_{w(n)}\rangle.$$ 
 The restriction map from the description in Theorem
\ref{th:cohFlag} to the fixed point set is as follows:
\begin{theorem}\label{th:corresp} Let $p_w\in Fl(V)$ be in the
fixed point set $Fl(V)^T$ as above. 
The inclusion $r_w:p_w\rightarrow Fl(V)$ induces a
restriction
\begin{equation} 
r^*_w: H^*_T(Fl(V)) \longrightarrow H^*_T(p_w) = 
S(\mathfrak{t}^*) =\C[u_1,\dots,u_n] 
\end{equation} 
such that 
$r_w^*: x_i \mapsto u_{w(i)}$ and $r^*_w: u_i \mapsto u_i$,
where $x_i$ and $u_i$, $i=1,\dots,n$ are the generators of the
equivariant cohomology in Theorem
(\ref{th:cohFlag}).  
In particular, $r:Fl(V)^T \hookrightarrow Fl(V)$ induces a map  
$$ 
r^*:H^*_T(Fl(V))\longrightarrow H^*_T(Fl(V)^T) = \bigoplus_{p\in W} 
\C[u_1,\dots, u_n] 
$$ 
whose further restriction to each component in the direct sum is 
$r_w^*$. 
\end{theorem}

\begin{demo}{Proof}
The classes $u_i$ come from the module structure of $H_T^*(M)$ over
$H_T^*$. For any $p\hookrightarrow M$, the induced map
$H_T^*(M)\rightarrow H_T^*(p)$ is a module map, and hence
$r_w^*u_i = u_i$
for all $w, i$.
A fixed point $p\in M$ 
in the Borel construction corresponds to a fixed copy of $BT\iso
p\times_T ET$ in
$M\times_TET$.   
We can restrict the $x_i$ to the 
fixed points of $Fl(V)$, and then use the Borel construction to make
the restriction equivariant.  

We first compute $r_w^*(x_1)$  for any 
$w\in W $, and proceed inductively to obtain $r_w^*(x_i)$ 
for all $i$.  Let $E:= V\times_T \C^n$ and note that $Fl(E) =
Fl(V)\times_TET$. The class $x_i=c_1(S_i)\in H^*(Fl(E))$ by definition. The
equivariant restriction to $p_w$ is the usual restriction to
$p_w\times_TET$ in $Fl(E)$. $S_1$ is the (pullback of) the
tautoligical line bundle over $\P(E)$. Under the projection to
$\P(E)$, $p_w\times_TET$ projects to $l_{w(1)}\times_TET$,
or a section $s_w$ of $\P(E)\rightarrow BT$. Then
$$
r_w^*(x_1)=c_1(S_1)|_{l_{w(1)}\times_TET}. 
$$
As a line bundle, $l_{w(1)}\times_TET = s_w^*S_1$ from the
commutative diagram
\begin{equation*}
\begin{CD}
s_w^*S_1 @>>> S_1\\
@VVV      @VVV\\
BT @>s_w>> \P(E).
\end{CD}
\end{equation*}
As $E$ splits into 
$$
E=\oplus_{i=1}^n L_i,
$$ with $L_i = l_i\times_{T_i}ET_i$  and $T_i \iso
S^1$ acts on $l_i$ with weight one and on $l_j$ trivially for $j\neq
i$, we have
$$ 
s_w^*S_1 = l_{w(1)}\times_TET=L_{w(1)}
$$ so that
$$
r_w^*(x_1)= c_1(L_{w(1)})
$$ by naturality of the pullback map.
We note that by definition, 
$ 
u_i=c_1(L_i)$ and hence 
$$ 
r_w^*(x_1)=c_1(S_1)|_{l_{w(1)\times_TET}}=c_1(L_{w(1)}) =u_{w(1)}. 
$$ 
We continue inductively: the point $p_w$ also specifies a 
two-dimensional fixed  
subspace $\langle l_{w(1)}, l_{w(2)}\rangle$ of $V$ 
containing $l_{w(1)}$. Under the Borel contruction, 
this two-space is a section in the projectivization $\P(Q_1)$ of the
quotient $Q_1=\pi_1^*E/S_1$, or the restriction of its tautological
line bundle $S_2$ to the copy of $BT$ that is the image of this
section. Then   
$$ 
r_w^*(x_2)=c_1(S_2)|_{l_{w(2)}\times_TET} = u_{w(2)}. 
$$ 
It is clear how to proceed from here. For each $x_i$, there is a 
 projection from $p_w$ to the corresponding $i$-space, which can
 (after crossing with $ET$ and modding out by $T$)
 be realized as a line in $S_i$. We obtain 
$$ 
r_w^*(x_i)=c_1(S_i)|_{l_{w(i)}\times_TET} = u_{w(i)} 
$$ 
which concludes the proof. \qed\end{demo}

An example of carrying out such a restriction may be useful to the
reader. 
\begin{example}\label{ex:restriction}
Let $n=3$, and $\alpha = (x_1-u_1)(x_1-u_3)$ be a $T$-equivariant
cohomology class on $\mathcal{O}_\lambda$, a generic coadjoint orbit
of $SU(3)$.  There are six fixed points of $T$ acting on
$\mathcal{O}_\lambda$, labeled by permutations in $S_3$. Let
$\alpha|_w$ indicate the restriction of $\alpha$ to $\lambda_w$. In one-line
notation, the restrictions are:
\begin{align*}
\alpha|_{[123]} = 0 & \hspace{.2in} \alpha|_{[213]} = (u_2-u_1)(u_2-u_3)\\
\alpha|_{[132]} = 0 & \hspace{.2in}\alpha|_{[231]} = (u_2-u_1)(u_2-u_3)\\
\alpha|_{[312]} = 0 & \hspace{.2in}\alpha|_{[321]} = 0.
\end{align*}
\end{example}

\subsection{Divided Difference Operators and Double Schubert
Polynomials} \label{se:schubertpolys}

The divided difference operators mentioned in Section \ref{se:intro}
were introduced by Bernstein, Gelfand and
Gelfand \cite{BGG} and independently by Demazure \cite{Demazure}. They
take polynomials to polynomials, and send to 0 anything symmetric in the
variables of the operator. For this reason, they
act on the (ordinary) cohomology of $G/T$, with $T$ the maximal
torus of a compact Lie group $G$.  These operators act on the 
$T$-equivariant cohomology of $G/T$ as well; by applying divided
difference operators to a certain equivariant cohomology class of degree
equal to the dimension of $G/T$, one can generate a linear basis for
$H_T^*(G/T)$. From a combinatorial perspective, this is equivalent in
the $G=SU(n)$ case to generating double Schubert polynomials by the
divided difference operators applied to a ``determinant
polynomial''. This was first carried out by Lascoux and
Sch\"utzenberger \cite{LS:schubert}, see also \cite{MacDonald}. In
this section we make explicit the connection between the combinatorics
and the equivariant cohomology.  By allowing the Weyl group to act on
the cohomology ring, we obtain {\em permuted} double Schubert
polynomials, which provide Weyl group many linear bases for
$H_T^*(G/T)$ as a module over $H_T^*$, as we will describe in Section
\ref{se:linearbasis}.

\begin{definition}
Let $f(x_i,x_{i+1})$ be a polynomial of variables $x_i$ and $x_{i+1}$
and possibly other variables. For each $1\leq i\leq n$ we define the
{\em divided difference operator} $\partial_i$ as follows:
$$
\partial_i f (x_i,x_{i+1})  = \frac{f(x_i,x_{i+1})-f(x_{i+1},x_i)}{x_i-x_{i+1}}.
$$ The resulting function $\partial_i f$ is also a polynomial.
\end{definition}

Let $W$ be the Weyl group for $SU(n)$, i.e. $W= N(T)/T$. Then $W$ is
isomorphic to $S_n$, the permutation group on $n$
letters. Furthermore, $W$ is generated by {\em simple transpositions}
$s_i$ which interchange $i$ and $i+1$. Any element $w\in W$ can thus
be written as a product $w=s_{i_1}s_{i_2}\cdots s_{i_l}$. Whenever $w$
is written with $l$ minimum, we call the expression a {\em reduced
word} for $w$. For any such reduced expression, we can define the operator
$$
\partial_{s_{i_1}s_{i_2}\cdots s_{i_l}} =
\partial_{i_1}\partial_{i_2}\cdots \partial_{i_l}.
$$
It turns out that the resulting operator is independent of the choice
of reduced word for $w$ (see \cite{MacDonald}). We can thus define the divided
difference operator associated to the element $w\in W$ as
$$
\partial_w = \partial_{i_1}\partial_{i_2}\cdots \partial_{i_l}
$$
for any reduced word expression $w=s_{i_1}s_{i_2}\cdots s_{i_l}$.

We now consider the following set of polynomials generated from one
polynomial by the divided difference operators $\partial _w$ and
permuted by the Weyl group.
\begin{definition}
The {\em determinant polynomial}
$\Delta\in\C[x_1,\dots,x_n,u_1,\dots,u_n]$ is defined to be
$$
\Delta (x,u) = \prod_{i<j}(x_i-u_j).
$$
The {\em permuted double Schubert polynomials} $\mathfrak{T}_w^\tau$
are defined by successive application of divided difference operators
to $\Delta$, as follows. The identity double Schubert polynomials are
$$
\mathfrak{T}_w^{id} := \partial_{w^{-1}}\Delta.
$$
The permuted double Schubert polynomials are defined from the identity
ones:
$$
\mathfrak{T}_w^\tau (x,u):= \mathfrak{T}_{\tau^{-1}w}^{id}(x,u_\tau)
$$
where $u_\tau$ indicates the permutation of the $u$ variables by
$\tau$. 
\end{definition}
\begin{example}\label{ex:class} We compute $\mathfrak{T}_w^\tau (x,u)$
for $n=3$, with $\tau = [213]$ and $w=[231] = s_1s_2$ in one-line
notation. We have $\Delta(x,u) = (x_1-u_2)(x_1-u_3)(x_2-u_3)$, and
$\Delta(x,u_\tau) = (x_1-u_1)(x_1-u_3)(x_2-u_3)$. Then
\begin{align*}
\mathfrak{T}_w^\tau (x,u) &= \mathfrak{T}_{[132]}^{[123]}(x,u_{[213]})\\ 
& =
\partial_2\Delta(x,u_{[213]})\\
& = \partial_2 (x_1-u_1)(x_1-u_3)(x_2-u_3)\\
&= (x_1-u_1)(x_1-u_3). 
\end{align*}
\end{example}

A couple comments are in order here.
\begin{enumerate}
\item For those familiar with double Schubert polynomials and determinant
polynomials, we note that in \cite{MacDonald} the determinant
polynomial is what we have called
$\Delta(x_{w_0}, u)$ for $w_0$ the long word in $W$. Double Schubert
polynomials $\mathfrak{S}_w(x,u)$ are equivalent to
$\mathfrak{T}_w^{w_0}$.
\item In this notation, the polynomials $\partial_w(x,u_\tau) =
\mathfrak{T}_{\tau w^{-1}}^\tau (x,u)$; equivalently,
$\mathfrak{T}_w^\tau = \partial_{w^{-1}\tau}(x,u_\tau)$.
\item While these definitions may seem to have an ungainly number of
inverses, the resulting geometric interpretation has a simple
statement; see Theorem \ref{th:support}.
\end{enumerate}

Under the quotient map
$$
\C[x_1,\dots,x_n, u_1,\dots,u_n]\longrightarrow H_T^*(G/T)
$$
the polynomials $\mathfrak{T}_v^\tau$ descend to cohomology classes; we
sloppily refer the classes themselves as $\mathfrak{T}_v^\tau$.

\subsubsection{The Bruhat Order and Permuted Schubert Varieties}
Let $G^\C$ be the complexification of a compact Lie group $G$, and let
$B\subset G^\C$ be a Borel subgroup. For example, for $G=SU(n)$,
$G^\C=Sl(n,\C)$ and one choice of $B$ is upper-triangular matrices. Then
$G^\C/B\iso G/T$.

The space $G^\C/B$ is composed of even-real-dimensional Schubert cells
indexed by elements in the Weyl group $W=N(T^\C)/T^\C$
$$
C_w:= BwB/B, w\in W.
$$
Technically, one needs to choose a lift of $w$ in $N(T^\C)$ but the
cell is independent of this choice, and so it is standard to consider
$w\in W$. The closures of these cells are {\em Schubert varieties}
$$
X_w := \overline{BwB/B}
$$
and were shown to generate the homology of $G^\C/B$ \cite{BGG}. We
define the {\em permuted Schubert cells} as
$$
C_w^\tau := \tau B\tau^{-1}wB/B
$$
and the {\em permuted Schubert varieties}
$$
X_w^\tau := \overline{\tau B\tau^{-1}wB/B}
$$
Note that these are $\tau B\tau^{-1}$-invariant varieties which
also generate the homology of $G^\C/B$.
\begin{example} The Schubert variety $X_{w_0}$ consists of all of
$G^{\C}/B$, whereas the variety $X_{id}$ consists of one point (the
identity coset). Similarly, $X_{\tau w_0}^\tau$ is all of $G^{\C}/B$ and  $X_\tau^\tau$ consists of just a point:
the permutation $\tau$.
\end{example}

The definition of these varieties suggest a partial ordering on
them, and hence on the elements of $W$ that index them. The Bruhat
order is defined:
$$
v\leq w \mbox{ if and only if } C_v\subset X_w.
$$
Similarly, we define a ``permuted'' ordering:
$$
v \leq_\tau w \mbox{ if and only if } C_v^\tau \subset X_w^\tau.
$$
The simple relation between these two is that
$$
v\leq_\tau w \mbox{ if and only if } \tau^{-1} v\leq \tau^{-1}w.
$$

Chevalley proved that the Bruhat ordering is equivalent to the
following.
\begin{definition}
We say that $v\leq w$ in the {\em Bruhat order} if and only if, for
all reduced word expressions $w=t_1\dots t_l$, there is a
subword $v=t_{i_1}\dots t_{i_k}$ with $i_1<\dots<i_k$ which is a
reduced word expression for $v$. We say that  $v\leq_\tau w$ in the
{\em permuted Bruhat order} if and only if $\tau^{-1} v\leq
\tau^{-1}w$. 
\end{definition}
Lastly, we note that the $T$ action on these varieties has fixed
points:
$$X_w^T = \{v\in W: v\leq w\}$$ and more generally,
$$(X_w^\tau)^T = \{v\in W: v\leq_\tau w\}.$$

We now show a fundamental relation between the permuted double
Schubert polynomials, viewed as cohomology classes, and the permuted
Schubert varieties.
\begin{theorem}\label{th:support}
The support of the equivariant cohomology class
$\mathfrak{T}_v^\tau$ is the permuted Schubert variety $X_v^\tau$.
\end{theorem}
\begin{example} This theorem states that 
\begin{align*}
supp\ \mathfrak{T}^{[213]}_{[231]} &= \{w\in W: w\leq_{[213]}
[231]\}\\
&= \{w\in W: [213]^{-1}w\leq [213]^{-1}[231]=[132]\} \\
&= \{w\in W: [213]^{-1}w=[123]\ or\ [132]\}\\
& = \{w = [213][123]=[213]\ or\ w=[213][132]= [231]\}.
\end{align*}
In Examples
\ref{ex:class} and \ref{ex:restriction}, where we showed that
$\mathfrak{T}^{[213]}_{[231]} = (x_1-u_1)(x_1-u_3)$ has support on
$\{[213], [231]\}$, as expected.
\end{example}
\begin{proof}
We first prove that $supp\ \mathfrak{T}_v^{id} = (X_v^{id})^T$ using Theorem \ref{th:corresp},
and then the theorem will follow easily from the definition of the
permuted double Schubert polynomials. Recall from Theorem
\ref{th:corresp} that the restrictions $u_i|_w = u_i$
and $x_i|_w = u_{w(i)}$ for all $w\in W$. We use
induction on $l(w)$. For $w=id = [1\ 2\ \cdots\ n]$, clearly the 
polynomial $\mathfrak{T}_{id}^{id}=\Delta$ restricts to zero at every
point except $\lambda_{id}$. Suppose now that the we have $supp\
\mathfrak{T}_w^{id} = (X_w^{id})^T$ for any $w$ with $l(w)=l-1$. Let
$l(v)=l$, and a reduced word expression for $v$ be
$v=s_{i_1}s_{i_2}\cdots s_{i_l}$. Let $w=vs_{i_l}$. Then restricted to
$\lambda_z$ we have
\begin{align*}
\mathfrak{T}^{id}_v|_{\lambda_z}&=\partial_{v^{-1}}\Delta|_{\lambda_z}=
\partial_{i_l}\cdots\partial_{i_1}\Delta |_{\lambda_z}=
\partial_{i_l}\mathfrak{T}_w^{id}|_{\lambda_z} \\
&=
\frac{\mathfrak{T}_w^{id}(x,u)-\mathfrak{T}_w^{id}(xs_{i_l},u)}{x_{i_l}-x_{i_{l+1}}}|_{\lambda_z}\\
&=
\frac{\mathfrak{T}_w^{id}(u_z,u)-\mathfrak{T}_w^{id}(u_{z
s_{i_l}} ,u)}{u_{z(i_l)}-u_{z(i_{l+1})}}.
\end{align*}
Suppose that $z\not\in (X_v^{id})^T$. Then $z\not\in(X_w^{id})^T$,
and since $l(w)=l-1$, by our inductive hypothesis the restriction
above is 
\begin{align*}
\mathfrak{T}^{id}_v|_{\lambda_z}&=-\frac{\mathfrak{T}_w^{id}(u_{z
s_{i_l}} ,u)}{u_{z(i_l)}-u_{z(i_{l+1})}}\\
&=-\frac{r^*_{z
s_{i_l}}\mathfrak{T}_w^{id}(x,u)}{u_{z(i_l)}-u_{z(i_{l+1})}}
\end{align*}
which is zero if $z s_{i_l}\not\in (X_w^{id})^T$.

Suppose then that $z s_{i_l}\in (X_w^{id})^T$. Then if $z<z
s_{i_l}$, we have $z\in (X_w^{id})^T$, a contradiction. If $z>z
s_{i_l}$, then $s_{i_l}$ increases the length of both $z s_{i_l}$
and $w$. But then $z s_{i_l}\in (X_w^{id})^T$ implies $z\in
(X_v^{id})^T$, again a contradiction.

We have shown that the restriction
$\mathfrak{T}^{id}_v|_{\lambda_z}$ is zero 
unless $z\in (X_v^{id})^T$. It is left to show that
$\mathfrak{T}^{\tau}_v|_{\lambda_z}=0$ unless $z\in
(X_v^{\tau})^T$. By definition,
$\mathfrak{T}^{\tau}_v(x,u):=
\mathfrak{T}_{\tau^{-1}v}^{id}(x,u_\tau)$, 
and we showed that $$supp\ \mathfrak{T}_{\tau^{-1}v}^{id}(x,u)=
(X_{\tau^{-1}v}^{id})^T = \{w\in W: w\leq \tau^{-1}v\}.$$
Then, permuting the $u$'s by $\tau$ we obtain
\begin{align*}
supp\ \mathfrak{T}^{\tau}_v &= supp\
\mathfrak{T}_{\tau^{-1}v}^{id}(x,u_\tau) \\
&= \{\tau w\in
W: w\leq \tau^{-1}v\}=\{w\in
W: \tau^{-1}w\leq \tau^{-1}v\}\\
& = (X_v^\tau)^T.
\end{align*}
\end{proof}

\subsection{A Linear Basis of $H_T^*(G/T)$}\label{se:linearbasis}
In \cite{Arabia:equivcohom} Arabia shows that there is a basis of
$H_T^*(G/T)$ as a module over $H_T^*$ with certain defining
properties. We use Theorem \ref{th:support} and Arabia's methods to
show that $\mathfrak{T}^\tau_w$ satisfy these properties for $G=SU(n)$. It follows
that, for any $\tau\in W$, $\alpha\in H_T^*(Fl(\C^n))$,
$$
\alpha = \sum_{w\in W} a_w^\tau \mathfrak{T}^\tau_w
$$
where $a_w^\tau\in H_T^*$.

\begin{theorem}\label{th:integral}
The classes $\mathfrak{T}_v^\tau$ have the (defining) properties
that
\begin{enumerate}
\item \ Their images in the regular cohomology of the flag variety
are a linear basis.
\item \ $\int_{X_w^{\tau w_0}}\mathfrak{T}_v^\tau=\delta_{wv}$,
where  $\int_{X_w^{\tau w_0}}$ is defined below.
\end{enumerate}
\end{theorem}
\begin{theorem}[Arabia]
A set of equivariant cohomology classes with properties (1) and (2) as
in Theorem \ref{th:integral} provides a linear basis for $H_T^*(G/T)$
as a module over $H_T^*$.
\end{theorem}

Following Arabia, the integral
$\int_{X_w^{id}}$ is defined by integrating over a certain (choice of)
smooth variety
which maps birationally to $X_w^{id}$.

\begin{theorem}[Bott-Samelson]
There is a (non-canonical) tower of $\C P^1$ bundles over $\C P^1$,
denoted $BS_{w_0}$, which maps birationally to $G/B$ such that, for
every $w\in W$, there is a smooth subvariety $BS_w$ of $BS_{w_0}$
which maps birationally to $X_w^{id}\subset G/B$.
\end{theorem}
While these {\em Bott-Samelson resolutions} are not canonical, we
always have that
\begin{equation*}
\begin{CD}
BS_w@>i_{BS_w}>> BS_{w_0}\\
@V\pi_wVV @VV\pi_{w_0}V\\
X_w @>i_w>> G/T
\end{CD}
\end{equation*}
is a $T$-equivariant commutative diagram (see
\cite{BS:G/Tring}), where $\pi_w$ is
generically one-to-one. In the induced commutative diagram
\begin{equation*}
\begin{CD}
H_T^*(BS_w) @<i^*_{BS_w}<< H_T^*(BS_{w_0})\\
@A\pi^*_wAA @AA\pi^*_{w_0}A\\
H_T^*(X_w) @<i^*_w<< H_T^*(G/T)
\end{CD}
\end{equation*}
the map $\pi^*_w$ is an injection for all $w$. Thus for any $\alpha\in
H_T^*(G/T)$, 
$$
\pi_w^*i_w^*\alpha = i_{BS_w}^*\pi_{w_0}^*\alpha.
$$
As each $BS_w$ is a smooth submanifold of $BS_{w_0}$, we define
$$
\int_{X_w^{id}}i_w^*\alpha := (\pi_w)_*\int_{BS_w}i^*_{BS_w}\pi^*_{w_0}\alpha
$$
where $(\pi_w)_*$ is the pushforward induced by $\pi_w$.
This whole construction can equally well be done for the varieties
$X_w^\tau$ to define $\int_{X_w^\tau}$ for all $\tau, w$.

\begin{demo}{Proof of Theorem \ref{th:integral}}
That the classes $\mathfrak{T}_w^\tau$ restrict to generators of the
regular cohomology of the flag manifolds (1) is by construction. The
$\mathfrak{T}_w^\tau$ restrict to Schubert polynomials (permuted by
$\tau$), shown in \cite{BGG} to be such generators.

The integrality condition (2) is also easy to show, using the
Atiyah-Bott/Berline-Vergne fixed point theorem.
Using Arabia's definition of the integral, the fixed
point theorem  says that
$$
\int_{X_w^{\tau w_0}} \mathfrak{T}_w^\tau = \sum_{p\in (BS_w^{\tau
w_0})^T} \frac{(\pi_{w_0}^*\mathfrak{T}^\tau_w)|_p}{e_{BS_w^{\tau w_0}}(p)}
$$
where $e_{BS_w^{\tau w_0}}(p)$ is the equivariant Euler class of $p\in
BS_w^{\tau w_0}$.
 The set of fixed
points for the $T$ action on $X_w^{\tau w_0}$ is
\begin{align*}
\{z: z\leq_{\tau w_0} w\} &= \{z: w_0\tau^{-1} z \leq
w_0\tau^{-1} w\}\\
&=\{z: \tau^{-1} z \geq\tau^{-1} w\}\\
&= \{z: z \geq_\tau w\}
\end{align*}
while $supp\ \mathfrak{T}_v^\tau = \{ z: z\leq_\tau v\}$. 

If $v=w$, there is only one point contributing to the integral, i.e.
$$
\int_{X_w^{\tau w_0}} \mathfrak{T}_w^\tau = \frac{\mathfrak{T}_w^\tau|_w}{e(w)}
$$
where $e(w)$ is the equivariant Euler class of $w\in BS_w^{\tau w_0}$
($BS_w^{\tau w_0}$ has exactly one point in the fibre over
$w\in G/T$). A quick computation shows that the
denominator of this expression is the product of the roots pointing
into the variety $X_w^{\tau w_0}$, while the numerator is the
product of the roots pointing out of $X_w^\tau$. As these are
equivalent, the quotient is 1.

If $v\neq w$, there are two possibilities. If $l(\tau^{-1}v)\leq l(\tau^{-1}w)$, then
the integral is clearly 0 because the sets $(X_w^{\tau w_0})^T$ and
$supp\ \mathfrak{T}_w^\tau$ do not intersect. If
$l(\tau^{-1}v)>l(\tau^{-1}w)$,  there may be a contribution to the
integral by points $z\in W$ such that $v\geq_\tau z \geq_\tau
w$. However, $\deg \mathfrak{T}_v^\tau = \deg
\mathfrak{T}_{\tau^{-1}v}^{id} = 
n(n-1)-l(\tau^{-1}v)$ and
\begin{align*}
\dim X_w^{\tau w_0} &= n(n-1)-\deg\mathfrak{T}_w^{\tau w_0} =
n(n-1)-[n(n-1) - l(w_0\tau^{-1}w)]\\
&= n(n-1)-l(\tau^{-1}w).
\end{align*}
Then $l(\tau^{-1}v)>l(\tau^{-1}w)$ implies $\dim  X_w^{\tau
w_0}>\deg \mathfrak{T}_v^\tau$, which implies that each contributing
term in the integral has a denominator of higher degree than the
numerator. Because the integral will be polynomial, these terms must
sum to zero.\qed
\end{demo}

\section{The Symplectic Picture of $H_T^*(\mathcal{O}_\lambda)$}
Recall from Theorem \ref{th:tolweis} that $M^\mu_\xi\subset M$ is the set
of points whose image under the moment map lies to one side of the
hyperplane $\xi^\perp_\mu$ through $\mu$ in $\mathfrak{t}^*$, i.e.
$$
M_\xi:= \{m\in M |\langle \phi(m),\xi\rangle \leq \langle\mu,\xi\rangle)\}.
$$ 

\begin{theorem}\label{th:supportononeside}  Let $\mathcal{O}_\lambda$
be a generic $SU(n)$ coadjoint orbit through $\lambda\in\mathfrak{t}^*$, and let
$\alpha\in H_T^*(\mathcal{O}_\lambda)$ 
be an equivariant cohomology class with $supp\ \alpha \subset
(\mathcal{O}_\lambda)^\mu_\xi$. Then there exists some $\tau\in W$ such
that
$$
\alpha=\sum_{v\in W} a_v^\tau \mathfrak{T}_v^\tau
$$ 
where $a_v^\tau\in H_T^*$ non-zero implies $supp\ 
\mathfrak{T}_v^\tau \subset (\mathcal{O}_\lambda)^\mu_\xi$. Furthermore,
$\tau$ may be chosen as any element of the Weyl group  such that $\xi$
realizes its minimum at $\phi(\lambda_\tau)$.
\end{theorem}

Note that this theorem is a symplectic, rather than topological
statement, i.e. it depends on the choice of $\lambda$ and $\mu$. We first need a lemma about the behavior of the $\xi$ with regard to
the Bruhat order.
We let
 $e_i$ be coordinate functions on $\mathfrak{t}^*$, so that if 
$\lambda\in\mathfrak{t}^*$ is written $(\lambda_1,\dots,\lambda_n)$, 
then $e_i(\lambda)=\lambda_i$. Recall that we have chosen $\lambda$ 
such that $\lambda_1>\dots> \lambda_n$.  
 Let $\xi^\perp_\mu\subset \mathfrak{t}^*$ indicate a  hyperplane
perpendicular to $\xi\in\mathfrak{t}$ through $\mu$. 

\begin{lemma}\label{le:order} Let $\xi\in \mathfrak{t}$ such that among points $\lambda_w$, $\xi$ attains its
minimum at $w=id$. Then $\xi$ respects the Bruhat order, i.e.
$\xi(\lambda_v)\leq\xi(\lambda_w)$ if $v\leq w$ in the Bruhat
order.
\end{lemma}
\begin{proof}
Write $\xi=\sum_{i=1}^n b_ie_i$ in the basis given above, with $b_i\in 
\mathbb{R}$. Comparing $\xi(\lambda_id)$ with $\xi(\lambda_{s_i})$, we have (by our minimality assumption)
$b_i\lambda_i+b_{i+1}\lambda_{i+1}\leq b_i\lambda_{i+1}+b_{i+1}\lambda_i$. Since
$\lambda_i>\lambda_{i+1}$, we find $b_i\leq b_{i+1}$. 
Over all $i$ we find
have
\begin{equation}\label{eq:b-order}
b_1\leq\cdots \leq b_n. 
\end{equation}
If $v<w$, then there is a sequence of length decreasing simple
reflections 
$s_{i_1}\cdots s_{i_k}$ such that $v=s_{i_1}\cdots s_{i_k}w$. For each 
reflection, we claim the value of $\xi$ decreases. Suppose $s_{i_k}=s_1$. 
Then 
$$\xi(\lambda_{s_1w})=\sum_ib_i\lambda_{w^{-1}s_1(i)} = 
b_1\lambda_{w^{-1}(2)}+b_2\lambda_{w^{-1}(1)}+
\sum_{i=3}b_i\lambda_{w^{-1}(i)}$$
and the difference $\xi(\lambda_w)-\xi(\lambda_{s_1w})$ is
$$
b_1\lambda_{w^{-1}(1)}+b_2\lambda_{w^{-1}(2)}-(b_1\lambda_{w^{-1}(2)}+b_2\lambda_{w^{-1}(1)}) 
= b_2(\lambda_{w^{-1}(2)}-\lambda_{w^{-1}(1)}) - b_1(\lambda_{w^{-1}(2)}-\lambda_{w^{-1}(1)}).
$$
But $s_1w<w$ if and only if $w^{-1}s_1<w^{-1}$, which implies that $w^{-1}(2)<w^{-1}(1)$, and thus
$\lambda_{w^{-1}(2)}>\lambda_{w^{-1}(1)}$. As the difference is positive, by
(\ref{eq:b-order}) we have
$$
b_1(\lambda_{w^{-1}(2)}-\lambda_{w^{-1}(1)})\leq b_2(\lambda_{w^{-1}(2)}-\lambda_{w^{-1}(1)})
$$
and therefore $\xi(\lambda_w)-\xi(\lambda_{s_1w})\geq 0$, as
desired.
The same proof applies for $s_{i_k}=s_j$ for any
length-reducing simple transposition $s_j$. Continuing
inductively we obtain the result.
\end{proof} 

This lemma generalizes quite easily to the case where the linear functional is minimized at $\lambda_\tau$ for any $\tau\in W$.

\begin{lemma}\label{le:sigmamin}
Let $\xi_\tau$ be a linear function on $\mathfrak{t}^*$ which attains its minimum on
$\lambda_\tau$. Then $v\leq_{\tau} w$ implies
$\xi_\tau(\lambda_v)\leq \xi_\tau(\lambda_w)$. 
\end{lemma}

\begin{proof} 
By definition, $v\leq_{\tau} w$ if and only if $\tau^{-1}v\leq
\tau^{-1}w$ in the Bruhat order. Then for any $\xi$ which is
minimized at $\lambda_{id}$, we have $\xi(\lambda_{\tau^{-1}v})\leq
\xi(\lambda_{\tau^{-1}w})$ by Lemma 
\ref{le:order}. Define 
$\xi (\lambda_w):= \xi_\tau(\lambda_{\tau w})$
 for all $w\in W$. Then
$\xi_\tau$ minimal at $\lambda_\tau$ implies $\xi$
minimal at $\lambda_{id}$, which then implies
$$\xi_\tau(\lambda_v)=\xi (\lambda_{\tau^{-1}v})\leq \xi(\lambda_{\tau^{-1}w})=\xi_\tau(\lambda_w).$$
\end{proof}

\begin{demo}{Proof of Theorem \ref{th:supportononeside}}
Let $\tau$ be such that  $\xi(\lambda_\tau)$ is minimal. For
any $\alpha\in H_T^*(M)$ we can write
$$
\alpha=\sum_{v\in W} a_v^\tau \mathfrak{T}_v^\tau
$$
with $a_v^\tau\in H_T^*$. We assume that $supp\ \alpha$ is contained
in $M_\xi^\mu$ and show that $supp\ \mathfrak{T}_v^\tau$ must also be
contained in $M_\xi^\mu$ whenever $a_v^\tau\neq 0$.

Suppose not. Let $F=\{ q\in W \mbox{such that } \alpha|_q = 0 \mbox{ but not all\
}a_v^\tau \mathfrak{T}_v^\tau|_q=0\}$. $F$ is not empty by assumption. Choose any $q\in F$
such that there are no points $q'\in F$ with $q'>_\tau q$. If
$a_q^\tau\mathfrak{T}_q^\tau=0$, then since $q\in F$, there exists
some $q'$ such that $a_{q'}^\tau\mathfrak{T}_{q'}^\tau|_q\neq
0$. Furthermore, $supp\ \mathfrak{T}_{q'}^\tau = (X_{q'}^\tau)^T$
implies $q\in (X_{q'}^\tau)^T$ and thus $q<_\tau q'$. Then
by Lemma \ref{le:sigmamin}, $\xi_\tau(q)\leq \xi_\tau(q')$, which
implies that $q'\in F$ since $supp\ \alpha\subset M_\xi^\mu$. Then $q$ was not a maximal element (in the
$>_\tau$ ordering) of $F$.
We now have that
$\alpha|_q = \sum_{v\in W} a_v^\tau \mathfrak{T}_v^\tau |_q=0$ implies 
$$
a_q^\tau\mathfrak{T}_q^\tau = -\sum_{v\neq
q}a_v^\tau\mathfrak{T}_v^\tau|_q \neq 0.
$$
By the same reasoning, for $v>_\tau q$, we have $a_v^\tau =0$ (otherwise $v\in F$ and
by Lemma \ref{le:sigmamin}, $q$ would not be maximal in the $>_\tau$
order). This implies 
$$
\sum_{v\ngeq_\tau
q}a_v^\tau\mathfrak{T}_v^\tau|_q \neq 0.
$$
But  $supp\ \mathfrak{T}_v^\tau = (X_v^\tau)^T$ implies $q\in
X_v^\tau$ for some $v$ for this sum to be non-zero, which in turn implies $q\leq_\tau v$, a
contradiction.\qed
\end{demo}

\section{Proof of Theorems 1.1 and 1.2}
First consider the case where $\lambda$ is generic (Theorem
\ref{th:generic}). 
We show that the kernel of the map
$$
\kappa_\mu:H_T^*(\mathcal{O}_\lambda) \longrightarrow
H^*(\mathcal{O}_\lambda/\!/T(\mu)) 
$$
is generated by the set of $\partial_v(x,u_\tau)$
listed in
Theorem \ref{th:generic}. The theorem then follows by the quotient
relation
$$
H^*(\mathcal{O}_\lambda/\!/T(\mu)) = H_T^*(\mathcal{O}_\lambda)/\ker
\kappa_{\mu}
$$
and the direct computation of $H_T^*(\mathcal{O}_\lambda)$ in Section
\ref{se:flagpresentation}.

We go about this by using support considerations. Consider the class
$\mathfrak{T}^\tau_w$ with support $(X_w^\tau) = \{v\in W:
v\leq_\tau w\}$. The functions 
$$\eta_k^\tau = \sum_{i=k+1}^n
e_{\tau(i)}$$
 obtain their maxima at $v=\lambda_w$ and their minima at
$v=\lambda_\tau$. Thus
$\eta_k^\tau(\lambda_w) = \sum_{i=k+1}^ne_{\tau(i)}(\lambda_w)<\sum_{i=k+1}^n\mu_{\tau(i)}$
implies $\sum_{i=k+1}^ne_{\tau(i)}(\lambda_v)<\sum_{i=k+1}^n\mu_{\tau(i)}$
by Lemma \ref{le:sigmamin}. We calculate
\begin{align*}
\sum_{i=k+1}^ne_{\tau(i)}(\lambda_w) &=
\sum_{i=k+1}^ne_{\tau(i)}((\lambda_{w^{-1}(1)},\dots,
\lambda_{w^{-1}(n)}))\\
&= \sum_{i=k+1}^n\lambda_{w^{-1}\tau(i)}.
\end{align*}
If $\sum_{i=k+1}^n\lambda_{w^{-1}\tau(i)}
 <\sum_{i=k+1}^n\mu_{\tau(i)}$, then 
 $supp\ \mathfrak{T}^\tau_w = 
(X_w^\tau)^T = \{v:v\leq_\tau w\}$ lies in
$M_{\eta_k^\tau}^\mu$, which by Theorem \ref{th:tolweis} proves that
$\mathfrak{T}^\tau_w$ is in $\ker \kappa_\mu$.
As $\mathfrak{T}^\tau_w (x,u)=
\partial_{w^{-1}\tau}\Delta(x,u_\tau)$, we have that $\partial_v\Delta(x,u_\tau)\in
\ker\kappa_\mu$ if $\sum_{i=k+1}^n \lambda_{v(i)} <
\sum_{i=k+1}^n\mu_{\tau(i)}$, as stated in Theorem \ref{th:generic}.

We need to show that these $\mathfrak{T}^\tau_w$ generate the kernel.
  Let
$\alpha\in H_T^*(M)$ be a homogeneous class such that $supp\ \alpha
\subset M_\xi^\mu$. By Theorem \ref{th:supportononeside} we can write
$$
\alpha = \sum_{w\in W} a_w^\tau\mathfrak{T}^\tau_w
$$
where $supp\ \mathfrak{T}_w^\tau\subset (X^\tau_w)^T=\{v:
v\leq_\tau w\}$.

We show $supp\ \mathfrak{T}_w^\tau\subset M_{\eta_k^\tau}^\mu$ for
some $k$.
Again by Lemma \ref{le:sigmamin}, if $v\leq_\tau
w$, we have $\eta^\tau_k(\lambda_v)\leq \eta^\tau_k(\lambda_w)$ for
all $k$. It is thus equivalent to show that
\begin{equation}\label{eq:eta_order}
\eta^\tau_k(\lambda_w)< \eta^\tau_k(\mu)
\end{equation} for some $k$.

Suppose that the equality (\ref{eq:eta_order}) does not hold for any
$k$. We have a series of inequalities 
\begin{align*}
\lambda_{w^{-1}\tau(n)} &\geq\mu_{\tau(n)}\\
\lambda_{w^{-1}\tau(n-1)} +
\lambda_{w^{-1}\tau(n)}&\geq\mu_{\tau(n-1)} +
\mu_{\tau(n)}\\ 
&\vdots\\
\lambda_{w^{-1}\tau(2)} + \dots + \lambda_{w^{-1}\tau(n)}&\geq\mu_{\tau(2)} + \dots + \mu_{\tau(n)}.
\end{align*}

Note that $supp\ \alpha\subset M_\xi^\mu$ is equivalent to
$$
\lambda_w\in supp\ \alpha \ \ \mbox{ implies }\ \ \xi(\lambda_w)<\xi(\mu).
$$
For $\xi=\sum_{i=1}^n b_ie_i$ by the same argument as that used in the
proof of Lemma \ref{le:order} we find that  $b_{\tau(1)}\leq\dots\leq
b_{\tau(n)}$. Therefore,  
\begin{align*}
(b_{\tau(n)}-b_{\tau(n-1)})\lambda_{w^{-1}\tau(n)}
&\geq (b_{\tau(n)}-b_{\tau(n-1)}) \mu_{\tau(n)}\\
(b_{\tau(n-1)}-b_{\tau(n-2)})(\lambda_{w^{-1}\tau(n-1)}
+ \lambda_{w^{-1}\tau(n)})&\geq (b_{\tau(n-1)}-b_{\tau(n-2)})(\mu_{\tau(n-1)} + \mu_{\tau(n)})\\
&\vdots\\
(b_{\tau(2)}-b_{\tau(1)})(\lambda_{w^{-1}\tau(2)} +
\dots + \lambda_{w^{-1}\tau(n)})& \geq
(b_{\tau(2)}-b_{\tau(1)}) (\mu_{\tau(2)} + \dots + \mu_{\tau(n)}).
\end{align*}
Summing the inequalities and using 
$\sum_{i=1}^n\lambda_i = \sum_{i=1}^n\mu_i=0$ we obtain
$$
b_{\tau(n)}\lambda_{w^{-1}\tau(n)} + \dots
+b_{\tau(1)}\lambda_{w^{-1}\tau(1)}\geq 
b_{\tau(n)}\mu_{\tau(n)}+\dots +
b_{\tau(1)}\mu_{\tau(1)},
$$ 
which is of course equivalent to
$$
b_n\lambda_{w^{-1}(n)}+\dots+ b_1\lambda_{w^{-1}(1)} \geq
b_n\mu_n+\dots + b_1\mu_1, 
$$
or $\xi(\lambda_w)\geq \xi(\mu)$, a
contradiction. 

It follows as an immediate corollary that the Tolman-Weitsman Theorem
(\ref{th:tolweis})
can be refined in the following way. Let $J\subset \mathfrak{t}$ be the set of
fundamental weights, permuted by the Weyl group. 

\begin{corollary}\label{co:sufficient} 
Let $M=\mathcal{O}_\lambda$ be the coadjoint orbit of $SU(n)$ through
a generic choice of $\lambda\in\mathfrak{t}^*$. Let $M^\mu_\xi$ and $K_\xi$ be as in
Theorem \ref{th:tolweis}. The kernel of the map
$$
\kappa_\mu:H_T^*(M)\longrightarrow H_T^*(M/\!/T(\mu))
$$ 
is the ideal generated by $K=\bigcup_{\xi\in J} K_\xi$. Equivalently, a
sufficient set of hyperplanes in Theorem \ref{th:tolweis} is the set of
all hyperplanes through $\mu\in\mathfrak{t}^*$ which are 
parallel to the codimension-one walls of the moment polytope.
\end{corollary}

\begin{proof} 
The fundamental weights of $SU(n)$ are $\eta_k^\tau$. We show how
these elements are perpendicular to hyperplanes parallel to
codimension-one walls in the image of the moment map.
For the moment map $\phi: M\rightarrow \mathfrak{t}^*$, 
codimension-one walls of the moment polytope consist of the image in
$\mathfrak{t}^*$ of the set of points fixed by some $S^1\subset T$ and having an
effective $T/S^1$ action.

The fixed point set of a codimension-one wall for
$\mathcal{O}_\lambda$ is the permutations of a partition of $n$
letters into two sets, of cardinality $k$ and $n-k$, respectively. One
easily sees that $\eta_k^\tau= \sum_{i=k+1}^n e_{\tau(i)}$ is
constant on the partition $(\{\tau(1),\dots,
\tau(k)\}, \{\tau(k+1),\dots, \tau(n)\})$ and its permutations
by $S_k\times S_{n-k}$. Thus $(\eta_k^\tau)^\perp$ is parallel to
this wall.
\end{proof}

We proceed to prove Theorem \ref{th:grassmannian}. We prove a little
lemma which shows that, as a subring of $H_T^*(Fl(\C^n))$, the
cohomology $H_T^*(Gr(k,n))$ is linearly generated by classes
$\mathfrak{T}^\tau_w$ that are symmetric in certain variables.
\begin{lemma}
Let $\alpha\in H_T^*(Gr(k,n))$. Then $\alpha$ can be written
$$
\alpha = \sum_{w\in W} a_w^\tau \mathfrak{T}_w^\tau
$$ with $\mathfrak{T}_w^\tau$ symmetric in
$(x_{\tau^{-1}(1)},\dots, x_{\tau^{-1}(k)})$ and in
$(x_{\tau^{-1}(k+1)},\dots, x_{\tau^{-1}(n)})$. 
\end{lemma}
\begin{proof} Let $i: Fl(\C^n)\rightarrow Gr(k,n)$ be the forgetful
map and $i^*: H_T^*(Gr(k,n))\rightarrow H_T^*(Fl(\C^n))$ be the
induced inclusion in equivariant cohomology.
The map $Fl(\C^n)^T\rightarrow Gr(k,n)^T$ under an identification
$Fl(\C^n)^T \iso W$ sends $S_n\rightarrow S_n/(S_k\times
S_{n-k})$. Any $\alpha\in H_T^*(Fl(\C^n))$ which is the image of $i^*$
is therefore constant in its restriction to a fixed point $p\in
Fl(\C^n)^T$ and its permutations by $S_k\times S_{n-k}$. By Section
\ref{se:linearbasis} we may write
$$
\alpha = \sum_{w\in W} a_w^\tau \mathfrak{T}_w^\tau.
$$ We need to show that $\mathfrak{T}_w^\tau$ are symmetric in the
relevant variables, which is equivalent to being in the image of
$i^*$. We show this first for $\tau = id$.

Suppose not all $\mathfrak{T}_w^{id}$ were symmetric in $x_1,\dots, x_k$
and in $x_{k+1},\dots, x_n$. Let $v\in W$ be the longest length
element such that $a_v^{id}\neq 0$ and $\mathfrak{T}_v^{id}$ not
symmetric in this sense. Then $v$ must be the longest element in $W$
in the orbit of $S_k\times S_{n-k}\cdot v$. If not, then since
$\alpha$ must be equal on all points in $S_k\times  S_{n-k}\cdot v$,
there must be some $\mathfrak{T}_w^{id}$ with $a_w^{id}\neq 0$,
$\mathfrak{T}_v^{id}\neq \mathfrak{T}_w^{id}$ and $l(w)>l(v)$. But
$\mathfrak{T}_w^{id}$ symmetric (as $v$ is the longest element with
$\mathfrak{T}_v^{id}$ not symmetric) implies $\alpha$ not symmetric.

For $v$ the longest element in the orbit, however,
$\mathfrak{T}_v^{id}$ is symmetric. For any $s_i\in S_k\times
S_{n-k}$, $s_iv<v$ implies $\partial_i\mathfrak{T}_v^{id}=0$, or
$\mathfrak{T}_v^{id}$ is symmetric in $x_i$ and $x_{i+1}$, $i=1,\dots,
k-1,k+1,\dots, n-1$.

Similarly, one proves for every $\tau$ that $\alpha = \sum_{w\in W}
a_w^\tau \mathfrak{T}_w^\tau$ implies that each contributing term
$\mathfrak{T}_w^\tau$ is symmetric in $(x_{\tau^{-1}(1)},\dots,
x_{\tau^{-1}(k)})$ and in 
$(x_{\tau^{-1}(k+1)},\dots, x_{\tau^{-1}(n)})$.
\end{proof}

It now follows that for $\alpha\in H_T^*(Gr(k,n)_\nu)$, Theorem
\ref{th:supportononeside} holds with all $\mathfrak{T}_v^\tau$
symmetric in the appropriate variables. Theorem \ref{th:grassmannian}
then follows by the same proof as that for Theorem \ref{th:generic}.

\section{Examples}
We present two examples for $n=4$, one of a generic coajoint orbit with two
positive and two negative eigenvalues, and the second of a degenerate
coadjoint orbit which is the 2-Grassmannian in $\C^4$.

\subsection{SU(4) Generic Coadjoint Orbits}
It has been shown that for $SU(n)$ coadjoint orbits that the reduction
at 0 of orbits whose isospectral sets $\langle
\lambda_1,\dots,\lambda_n\rangle$ are close to $\langle
a,\dots,a,b\rangle$ or $\langle a,b,\dots,b\rangle$ for $a>b$ are
again coadjoint orbits, now of $SU(n-1)$\cite{GLS:fibrations}.
For $0$-weight varieties, this means that if $n-1$ eigenvalues are
above (or below) zero, the associated $0$-weight variety will be a
coadjoint orbit of $SU(n-1)$. 
Here we use the
methods developed above to compute
the cohomology ring
of the $0$-weight variety of $SU(4)$ in
the case where $\lambda_1>\lambda_2>0>\lambda_3>\lambda_4$.

There is one choice that has an effect on the kernel of the Kirwan map
$\kappa:H^*_T(\mathcal{O}_\lambda)\longrightarrow
H^*(\mathcal{O}_\lambda/\!/T(0))$: either $\lambda_2+\lambda_3>0$
(which implies $\lambda_1+\lambda_4<0$), or  $\lambda_2+\lambda_3<0$
(which implies $\lambda_1+\lambda_4>0$). In this case, the resulting cohomology
rings are isomorphic, and so we choose $\lambda_2+\lambda_3>0$. This
forces the ordering on the partial sums:
\begin{equation}\label{eq:partialsum}
\lambda_1+\lambda_2>\lambda_1+\lambda_3>\lambda_2+\lambda_3>0>\lambda_1+\lambda_4>\lambda_2+\lambda_4>\lambda_3+\lambda_4.
\end{equation}

\begin{theorem}\label{th:nonextremal}
For $\mathcal{O}_\lambda$ the coadjoint orbit through $\lambda$
satisfying the relation (\ref{eq:partialsum}), the cohomology
of the zero weight variety is described by
$$
H^*(\mathcal{O}_\lambda/\!/T(0))= \frac{\mathbb{C}[x_1,\dots,x_4,u_1,\dots,u_4]}{\bigg(\prod(1+u_i)-\prod(1+x_i),\sum_iu_i,\alpha_1,\dots,\alpha_{14}\bigg)}
$$ where $\alpha_i$ for $i=1,\dots,14$ are the following degree 4 classes:
\begin{align*}
\alpha_1 &=(x_1-u_4)(x_2-u_4)\hspace{.3in} \alpha_3 =(x_1-u_2)(x_2-u_2)\\
\alpha_2 &=(x_1-u_3)(x_2-u_3)\hspace{.3in} \alpha_4 =(x_1-u_1)(x_2-u_1)\\ 
\alpha_{5}&= x_1^2+x_1x_2+x_2^2 - (u_2+u_3+u_4)(x_1+x_2) +
(u_2u_3+u_2u_4+u_3u_4)\\ 
\alpha_{6} &=x_1^2+x_1x_2+x_2^2 - (u_1+u_3+u_4)(x_1+x_2) +
(u_1u_3+u_1u_4+u_3u_4)\\ 
\alpha_{7} &=x_1^2+x_1x_2+x_2^2 - (u_1+u_2+u_4)(x_1+x_2) +
(u_1u_2+u_1u_4+u_2u_4)\\ 
\alpha_{8} &=x_1^2+x_1x_2+x_2^2 - (u_1+u_2+u_3)(x_1+x_2) + (u_1u_2+u_2u_3+u_1u_3)\\
\alpha_{9} &= x_1x_2+x_2x_3+x_1x_3 - (u_3+u_4)(x_1+x_2+x_3) + (u_3^2 +
u_3u_4+ u_4^2)\\
\alpha_{10} &=x_1x_2+x_2x_3+x_1x_3 - (u_2+u_4)(x_1+x_2+x_3) + (u_2^2 +
u_2u_4+ u_4^2)\\
\alpha_{11} &=x_1x_2+x_2x_3+x_1x_3 - (u_1+u_4)(x_1+x_2+x_3) + (u_1^2 +
u_1u_4+ u_4^2)\\
\alpha_{12} &=x_1x_2+x_2x_3+x_1x_3 - (u_2+u_3)(x_1+x_2+x_3) + (u_2^2 +
u_2u_3+ u_3^2)\\
\alpha_{13} &= x_1x_2+x_2x_3+x_1x_3 - (u_1+u_3)(x_1+x_2+x_3) + (u_1^2 +
u_1u_3+ u_3^2)\\
\alpha_{14} &=x_1x_2+x_2x_3+x_1x_3 - (u_1+u_2)(x_1+x_2+x_3) + (u_1^2 +
u_1u_2+ u_2^2)\
\end{align*}
\end{theorem}
\begin{corollary}
Let $\mathcal{O}_\lambda$ be a non-extremal coadjoint orbit of
$SU(4)$. The Poincar\'e polynomial for $\mathcal{O}_\lambda/\!/T(0)$
is $1+6t^2+6t^4+t^6$.
\end{corollary}

For this very symmetric case, one can check the Betti numbers
using the Morse theory methods of Kirwan; see
 \cite{KiBook:quotients}, where she does similar examples. For the
Betti numbers of this example, it is easier to use Kirwan's methods,
because the indices
 of critical points are 
 easy to calculate using the symmetry.   When the
 entries in $\lambda$ are distributed randomly, however,  it may be quite
difficult 
 to read off the index of a particular critical point. The method
 proposed here, while requiring many more tedious calculations, has the
 advantage of being straightforward.

\begin{demo}{Proof of Theorem}
First consider $\tau=id$. The fundamental weights which are
minimized at $\lambda_{id}$ (among $\lambda_v$) are $\xi=e_4,
\xi=e_3+e_4$ and $\xi= e_2+e_3+e_4$. The set of points $\lambda_w$
such that $\xi(\lambda_w)<\xi(0)$ for one of these choices of $\xi$
are as follows, listed by length:

\begin{gather*}
\begin{aligned}
\mbox{For }l(w)=0, w&=[1234]\\
\mbox{For }l(w)=1, w&=[2134], [1243], [1324]\\
\mbox{For }l(w)=2, w&=[2143], [1423], [3124], [2314], [1342]\\
\mbox{For }l(w)=3, w&=[3214], [3142],  [2341], [2413], [4123]\\
\mbox{For }l(w)=4, w&=[3241], [4132], [4213].
\end{aligned}
\end{gather*}

Correspondingly, we need to find the set of all
$\mathfrak{T}_{w}^{\tau}= \mathfrak{T}_{w}^{id}$ for
these $w$. The smallest degree $\mathfrak{T}_{w}^{id}$ are those
for which $l(w)=4$. We compute these classes:

\noindent For $w=[3241]=s_1s_2s_1s_3$, 
\begin{gather*}
\begin{aligned}
\mathfrak{T}_{w}^{id}=\partial_{w^{-1}} \Delta&=
\partial_3\partial_1\partial_2\partial_1
(x_1-u_2)(x_1-u_3)(x_1-u_4)(x_2-u_3)(x_2-u_4)(x_3-u_4) \\
&= \partial_3\partial_1\partial_2
(x_1-u_3)(x_1-u_4)(x_2-u_3)(x_2-u_4)(x_3-u_4) \\
&= \partial_3\partial_1 (x_1-u_3)(x_1-u_4)(x_2-u_4)(x_3-u_4) \\
&= \partial_3(x_1-u_4)(x_2-u_4)(x_3-u_4)\\
&= (x_1-u_4)(x_2-u_4).
\end{aligned}
\end{gather*}
\noindent For $w=[4132]=s_3s_2s_3s_1$,
\begin{gather*}
\begin{aligned}
\mathfrak{T}_{w}^{id}=\partial_{w^{-1}} \Delta&=
\partial_1\partial_3\partial_2\partial_3
(x_1-u_2)(x_1-u_3)(x_1-u_4)(x_2-u_3)(x_2-u_4)(x_3-u_4) \\
&= \partial_1\partial_3\partial_2
(x_1-u_2)(x_1-u_3)(x_1-u_4)(x_2-u_3)(x_2-u_4) \\
&= \partial_1\partial_3 (x_1-u_2)(x_1-u_3)(x_1-u_4)(x_2+x_3-(u_3+u_4)) \\
&= \partial_1(x_1-u_2)(x_1-u_3)(x_1-u_4)\\
&= x_1^2+x_1x_2+x_2^2-(u_2+ u_3+u_4)(x_1+x_2) + (u_2u_3+u_2u_4+u_3u_4).
\end{aligned}
\end{gather*}
\noindent For $w=[4213]=s_3s_1s_2s_1$,
\begin{gather*}
\begin{aligned}
\mathfrak{T}_{w}^{id}=\partial_{w^{-1}} \Delta&=
\partial_1\partial_2\partial_1\partial_3
(x_1-u_2)(x_1-u_3)(x_1-u_4)(x_2-u_3)(x_2-u_4)(x_3-u_4) \\
&= \partial_1\partial_2\partial_1
(x_1-u_2)(x_1-u_3)(x_1-u_4)(x_2-u_3)(x_2-u_4) \\
&= \partial_1\partial_2 (x_1-u_3)(x_1-u_4)(x_2-u_3)(x_2-u_4) \\
&= \partial_1(x_1-u_3)(x_1-u_4)(x_2+ x_3-(u_3+u_4))\\
&= x_1x_2+x_1x_3+x_2x_3-(u_3+u_4)(x_1+x_2+x_3) + (u_3^2+u_3u_4+u_4^2).
\end{aligned}
\end{gather*}

These three degree 4 cohomology classes can be permuted by the Weyl
group to obtain other classes in the kernel, for if $\xi$ is minimized
at $\lambda_{id}, \xi=\sum_{i=1}^nb_ie_i,$ then $\xi_\tau =
\sum_{i=1}^nb_ie_{\tau(i)}$ is minimized at
$\lambda_\tau$. If $\xi(w)<\xi(0)=0$, then
$\xi_\tau(\lambda_{\tau w})=\xi(\lambda_w)<\xi(0)=\xi_\tau(0)=0$
implies that $\mathfrak{T}_{\tau w}^{\tau}\in \ker
\kappa$. But $\mathfrak{T}_{\tau w}^{\tau}:=
\mathfrak{T}_{w}^{id}(x,u_{\tau})$, or permutations in the
$u$-variables of the three classes we just computed. The permutations
of $\mathfrak{T}_{w}^{id}$ for $w=[3241]$ are the classes
$\alpha_1,\dots,\alpha_4$ listed in Theorem \ref{th:nonextremal}. The permutations
of $\mathfrak{T}_{w}^{id}$ for $w=[4132]$ are the classes
$\alpha_5,\dots, \alpha_{8}$, and those for $w=[4213]$ are the classes
$\alpha_{9},\dots,\alpha_{14}.$

These classes are independent contributions to the kernel of 
$$
\kappa:H_T^*(\mathcal{O}_\lambda)\longrightarrow
H^*(\mathcal{O}_\lambda/\!/T(0)),
$$
as can be verified by a laborious computation. While these classes
are the only degree 4 classes, we theoretically must calculate the
remaining (higher) degree classes in the kernel. However, a direct
computation shows that the Poincar\'e polynomial of
$H^*_T(\mathcal{O}_\lambda)/\langle \alpha_1,\dots,\alpha_{14}\rangle$
is $1+6t^2+6t^4+t^6$, and since the
$\dim_{\R}\mathcal{O}_\lambda/\!/T(0)=6$, there cannot be any further
contributions to the kernel. \qed\end{demo}
\begin{figure}
\centerline{
\psfig{figure=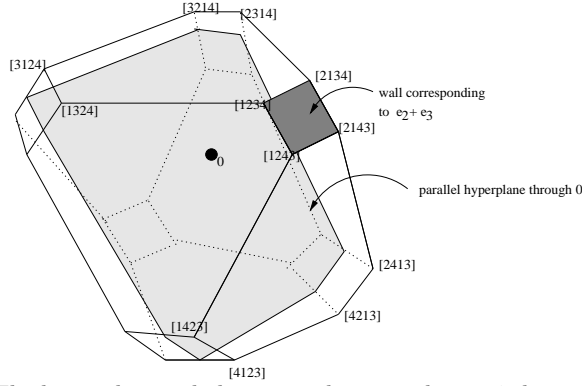,width=3.0in}
}
\centerline{
\parbox{5in}{\caption[A Hyperplane Parallel to a Wall for
$SU(4)$]{{\small The hyperplane in light gray indicates a degree 4
class in the kernel of the Kirwan map
$\kappa:H_T^*(\mathcal{O}_\lambda)\rightarrow
H^*(\mathcal{O}_\lambda/\!/T(0))$ whose support lies on the labeled
points to one side of the hyperplane.}}}
}
\end{figure}

\subsection{The Grassmannian $Gr(2,4)$}
The Grassmannian case is very similar in computation to the case of
the generic coadjoint orbit of $SU(4)$ with the eigenvalues ordered
above. As a (degenerate) coajoint orbit, $Gr(2,4)_\nu$ is the orbit
through a point $\nu\in\mathfrak{t}^*$ where the orbit has two distinct
eigenvalues instead of four. We label them
$\nu=(\nu_1,\nu_1,\nu_2,\nu_2)$, so that $\nu_1+\nu_1+\nu_2+\nu_2=0$,
or $\nu_1=-\nu_2$. This case is in some sense a limit case of the
example above, as $\lambda_1$ approaches $\lambda_2$, and $\lambda_3$
approaches $\lambda_4$ in Expression (\ref{eq:partialsum}).

We compute the cohomology of the symplectic reduction
$Gr(2,4)_\nu/\!/T(\mu)$ where the point of reduction $\mu\in\mathfrak{t}^*$
cannot be zero, as zero is not a regular value of the moment map in
this case. Choose $\mu=(\mu_1,\mu_2,\mu_3,\mu_4)$ with
$\sum_{i=1}^4\mu_i=0$ and $\mu_1>\mu_2>\mu_3>0$ but close enough to
zero that $\mu$ is in the image of the moment polytope.

\begin{theorem}
For $\mathcal{O}_\nu$ the coadjoint orbit through $\nu$
satisfying the relation above, the cohomology
of the $\mu$ weight variety for $\mu$ as above is described by
$$
H^*(Gr(2,4)_\nu/\!/T(\mu))= \frac{\mathbb{C}[x_1+x_2,x_3+x_4, x_1x_2,x_3x_4,u_1,\dots,u_4]}{\bigg(\prod(1+u_i)-\prod(1+x_i),\sum_iu_i,\beta_1,\beta_2,\beta_{3},\alpha_1,\dots,\alpha_8\bigg)}
$$ where the $\alpha_i$  are the degree 4 classes listed in Theorem
\ref{th:nonextremal} and $\beta_i$ are the following degree 2 classes:
\begin{align*}
\beta_1 &=(x_1+x_2-(u_1+u_2))\\
\beta_2 &=(x_1+x_2-(u_1+ u_4))\\
\beta_3 &=(x_1+x_2-(u_1+u_3))\\
\end{align*}
\end{theorem}
\begin{proof} It should first be noted that indeed there is repetition
in the classes $\alpha_i$ and $\beta_i$. This is done however to
emphasize that the classes in the kernel include the classes in the
kernel for the generic case which are in the image of the map
$H_T^*(Gr(2,4))\hookrightarrow  H_T^*(SU(4)/T)$ induced by the
forgetful map from the flag variety to the Grassmannian.

We do a similar computation as above to find classes in the kernel of
the map
\begin{equation}\label{map:grass}
\kappa_\mu: H_T^*(Gr(2,4)_\nu)\longrightarrow
H^*(Gr(2,4)_\nu/\!/T(\mu)).
\end{equation}

The eight classes
$\alpha_1,\dots,\alpha_8$ are easily
seen to be in the kernel of $\kappa_\mu$. Theorem \ref{th:generic}
states that if $\xi_\tau^j(\lambda_w)< \xi_\tau^j(0)$, then
$\mathfrak{T}_w^\tau\in \kappa_0$ for the case of a generic
coadjoint orbit reduced at 0.  But for $\mu$ small enough,
$\xi_\tau^j(\lambda_w)< \xi_\tau^j(0)$ implies
$\xi_\tau^j(\nu_w)<\xi_\tau^j(\mu)$, which implies that the same
cohomology class $\mathfrak{T}^\tau_w$ will be in the kernel of
$\kappa_\mu$ for the map (\ref{map:grass}) if it is an element of the
cohomology $H_T^*(Gr(2,4))$. The classes
$\alpha_1,\dots,\alpha_8$ are exactly
those classes in the kernel of $\kappa_0$ which are in the image of the injection
$H_T^*(Gr(2,4))\hookrightarrow  H_T^*(SU(4)/T)$.

To compute the remaining degree 2 classes, we do as before: 
Consider $\tau = [4321]$. A fundamental weight minimized at
$\nu_{[2211]}$ is $\xi=e_1+e_2$. We note that
$$
\xi(\nu_{[2211]})<
\xi(\nu_{[1212]})=\xi(\nu_{[2112]})=\xi(\nu_{[2121]})=\xi(\nu_{[1212]})=0
$$
  whereas $\xi(\mu)=\mu_1+\mu_2>0$. Note that the fixed points
  $\lambda_w$ in $\mathcal{O}_\lambda$
  which map to these points in $\mathcal{O}_\nu$ under the forgetful
  map do not have the same property as above. When one evaluates $\xi$
  on the points $\lambda_w$, several will be above $\xi(0)$, which is
  why we anticipate that the corresponding $\mathfrak{T}^\tau_w$
  will not have been already seen in the quotient for the generic
  case.

Let $w=[1324]$, which maps to $[1212]$ under the forgetful map, and
calculate $\mathfrak{T}^\tau_w (x,u) =
\mathfrak{T}_{\tau^{-1}w}^{id}(x,u_\tau)$. First we find
$\tau^{-1}w = [4231]=s_3s_1s_2s_1s_3$:
\begin{align*} 
\mathfrak{T}^{id}_{\tau^{-1}w} (x,u) &=
\partial_{w^{-1}\tau}\Delta\\ 
&= \partial_3\partial_1\partial_2\partial_1\partial_3\Delta\\
&= \partial_3\partial_1\partial_2\partial_1
(x_1-u_2)(x_1-u_3)(x_1-u_4)(x_2-u_3)(x_2-u_4)\\
&= \partial_3\partial_1\partial_2
(x_1-u_3)(x_1-u_4)(x_2-u_3)(x_2-u_4)\\
&= \partial_3\partial_1(x_1-u_3)(x_1-u_4)(x_2+x_3-(u_3+u_4))\\
&= \partial_3(x_1x_2+x_2x_3+x_1x_3 -
(u_3+u_4)(x_1+x_2+x_3)+u_3^2+u_3u_4+u_4^2)\\
&= (x_1+x_2)-(u_3+u_4).
\end{align*}

Then $\mathfrak{T}^\tau_w (x,u)= (x_1+x_2)-(u_2+u_1)=\beta_1$.
Note that $\mathfrak{T}^\tau_w (x,u)$ is symmetric in $x_1$ and
$x_2$, which is necessary for it to be an element of the
kernel.  Similarly, one finds that there are no other degree two
classes given by this choice of $\xi$.

For $\xi = e_2+e_3$ and $\xi=e_1+e_3$ and using the same techniques,
one finds the classes $\beta_2 = x_1+x_2-(u_1+u_4)$ and $\beta_3=
x_1+x_2-(u_1+u_3)$, respectively.

One can check that the quotient ring is actually the cohomology of the
two-sphere. In fact, it is known that $Gr(2,4)//T$ is always a
two-sphere for regular values of the moment map in the moment
polytope.

\end{proof}









\end{article}
\end{document}